# Archimedean cuboctahedron: The Medieval Journey from the Middle East to Northern Russia

**Andrey Yu. Chernov[1], Eugene A. Katz[2]**

[1] e-mail: trezin@yandex.ru;
[2] Swiss Institute for Dryland Environmental and Energy Research,, J. Blaustein Institutes for Desert Research, Ben-Gurion University of the Negev, 8499000, Sede Boqer Campus, Israel; e-mail: keugene@bgumail.bgu.ac.il

**Abstract**. Bronze *cuboctahedral* weights dated to the VIII-X centuries were found in northwest Russia near Ladoga, one of the most important trading centers in Eastern Europe in the VIII-X centuries. The history of the mathematical study of *cuboctahedron* and more generally of the entire family of *Archimedean solids* in the Middle East and Europe supports the archeological hypothesis about the origin of these artifacts and trading contacts between Europe and the Islamic Caliphate at that time when European mathematicians were not aware of such polyhedra, but Arab-Persian scientists and craftsmen were.

"The Mathematical Intelligencer" has repeatedly published articles on *Archimedean solids* [1-5]. In this essay, we comprehensively describe the complicated history of the geometrical exploration of these semiregular polyhedra to support an archeological hypothesis on the origin of bronze weight dated to the VIII-X centuries in the form of a cuboctahedron. The latter is one of the thirteen Archimedean polyhedra with 8 triangular and 6 square faces (Figure 1). It was found by the expedition of the Russian archaeologist Evgeny Ryabinin (1948 - 2010) working from 1973 to 2001 on the banks of Lake Ladoga and the Volkhov River (north-west of Russia) near the village of *Staraya Ladoga*, known as the town of *Ladoga* until 1704 (Figure 2). One of the authors (A Yu Ch) participated in this research.

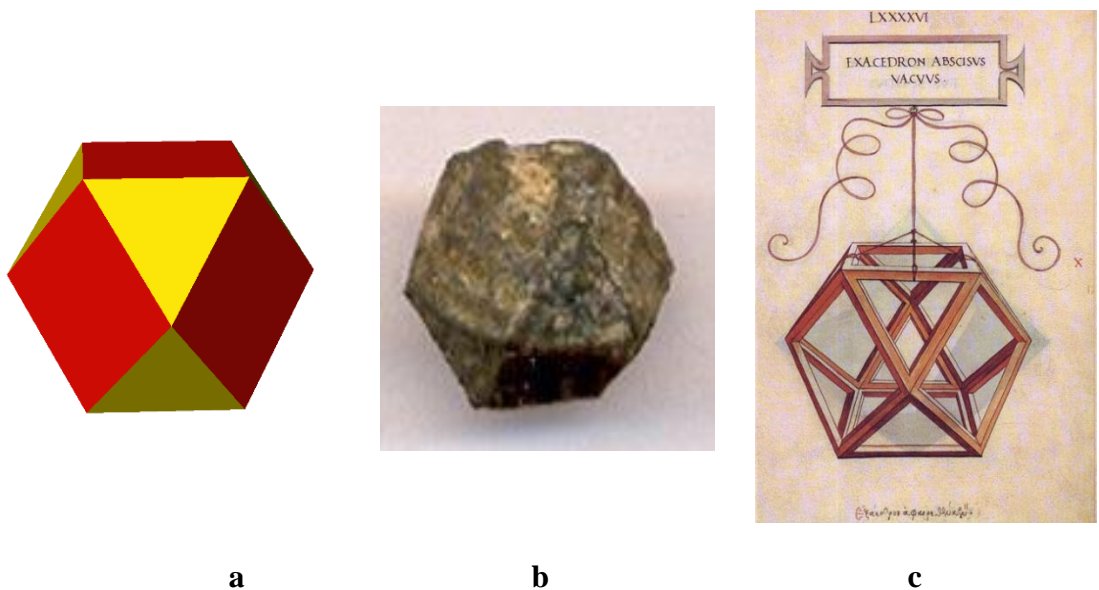

        **a**                 **b**                 **c**

Figure. 1. (a) Cuboctahedron. (b) The bronze weight of the VIII-X centuries found on the banks of Lake Ladoga and the Volkhov River. Photo by A. Yu. Chernov. (c) Leonardo da Vinci's drawing of the cuboctahedron ('*exacedron abscisus vacuus*') for Luca Pacioli's book "De Divina Proportione".

Ladoga is known as the oldest town in Russia and is distinguished by the unique preservation of an almost four-meter cultural layer [6]. The archaeological expedition by Evgeny Ryabinin established the dendrochronological date of the emergence of this village in the lower reaches of the Volkhov River as 753 [7-8].

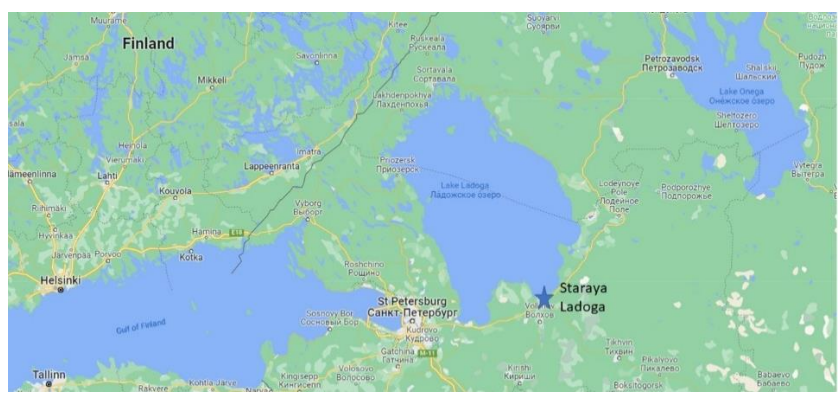

**a**

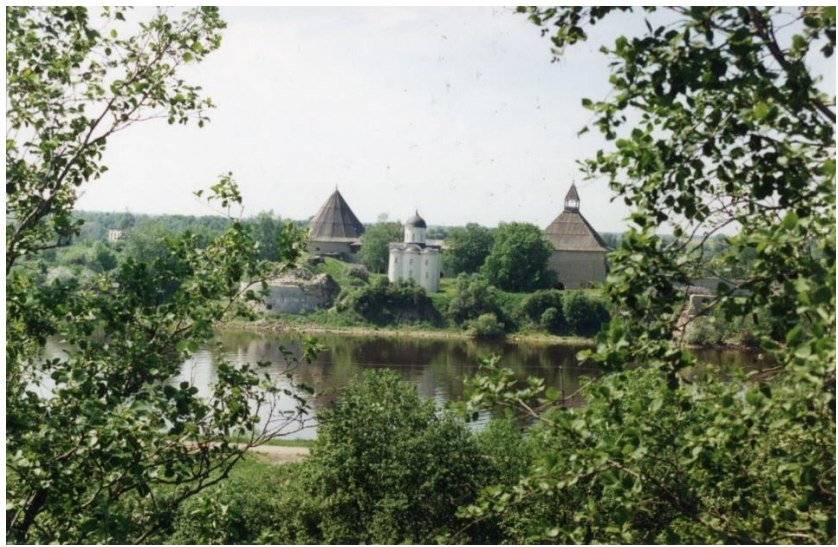

b

Figure 2. Staraya Ladoga. (a) geographical location; (b) Photo by A. Yu Chernov.

According to the *Hypatian Codex* (also known as *Hypatian Letopis* or *Ipatiev Letopis*; Russian: Ипатьевская летопись) [9], the semi-legendary Varangian leader Rurik, identified today with the Danish prince Rorik of Jutland and Friesland (810-880), arrived at Ladoga in 862. Rurik later moved to Novgorod and subsequently, his successors moved from there to Kiev where foundations for the powerful state of *Kievan Rus'* were laid [10]. Probably for that reason, Ladoga is considered sometimes the first capital of Russia.

Ladoga was one of the most important trading centers in Eastern Europe in the VIII-X centuries, and it is estimated that between 90% to 95% of all Arab dirhams found in Sweden passed through Ladoga. Merchant vessels sailed from the Baltic Sea through Ladoga to Novgorod and then to Constantinople. This route is known as the *trade route from the Varangians to the Greeks*. An alternative way led down the Volga River along the Volga trade route to the Khazar capital of Atil, and then to the southern shores of the Caspian Sea, all the way to Baghdad.

The remains of a blacksmith and jewelry workshop and a production complex, discovered by Ryabinin's expedition, are dated to the VIII – X centuries, the owners of which were engaged in the manufacture of the oldest glass jewelry in Eastern and Northern Europe. The factory used Arabic low-temperature technology to produce beads from the local quartz sand but with the addition of ash from salt marsh plants, that was brought from the Middle East. The Ladoga glass workshop operated from the 770s to the end of the 830s. Especially popular among the first settlers of Ladoga were the so-called *eye beads* - beautiful products with multi-colored spots - "*eyes*" (Figure 3).

There is a note about such eye beads in *the Tale of Bygone Years* (also known as *the Russian Primary Chronicle* [10-12]) which author visited Ladoga in 1114 and carried out the first archaeological research on this site. He was struck by the sight of such beads, which had fallen into disuse at least a century before him.

The financial and production scheme that the enterprising Ladoga residents developed and implemented was simple but effective: glass beads were transported by merchants in large birch bark boxes to the Finnish lands of the forested North, where they were exchanged with Finnish hunters for furs. The highly valued furs were then sold by the Ladoga residents to Arabic merchants from the Middle East in exchange for silver coins – *dirhams*: "Fur for beads. Coins for fur" (Figure 4). The glassworks existed in the mode of an emission bank: no matter how much furs entered the market, so many glass "eyes" were produced [13-14].

The eyed beads and medieval silver coins, including Arabic dirhems from the VIII-X centuries, can still be found on the banks of the Volkhov River.

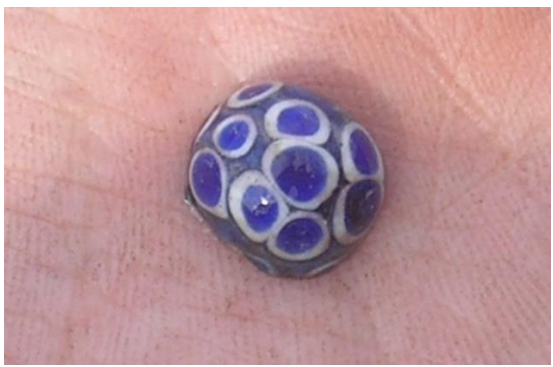
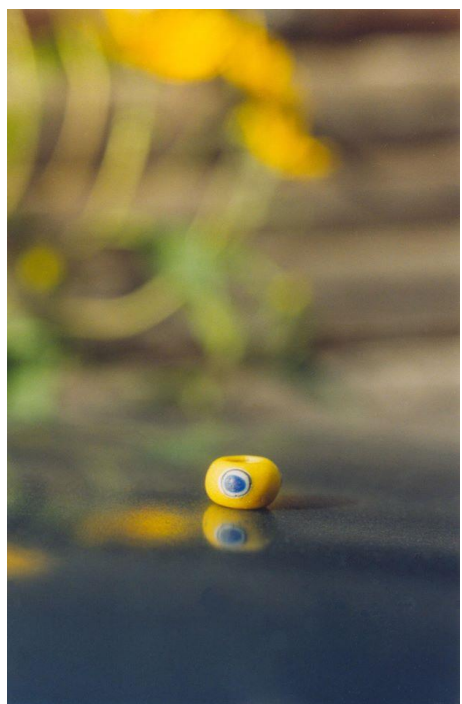

Figure 3. Medieval eye beads. Found in Old Ladoga and photographed by A. I. Chernov.

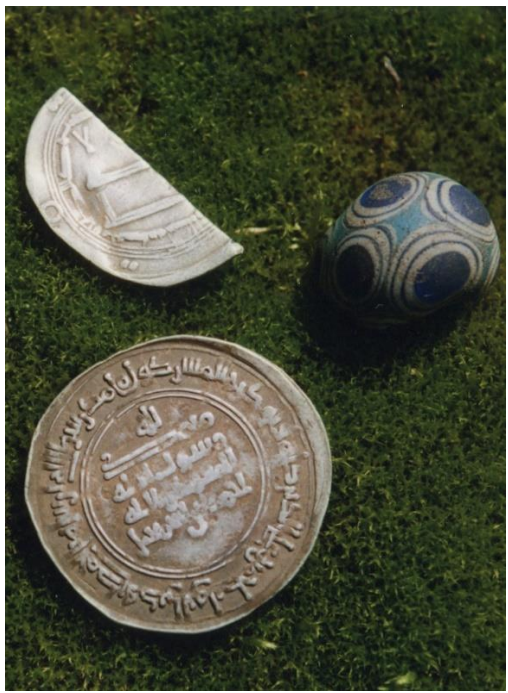

Figure 4. International trade currency of the VIII-IX centuries: Ladoga glass beads, Arabic dirhams, and fragmented dirhams. Photo by A. Yu. Chernov.

Where there is a developed exchange of commodities, there must be a measure of this exchange. Arab merchants had to determine the weight of their dirhams. The polyhedral weights found by archeologists in Ladoga along with dirhams (Figure 4) can be considered as clear evidence of the trade. One can see the reason for the need to measure the dirham weight at the point of the sale in Figure 4: the coins were simply fragmented by weight. Of course, the scales and weights were subsequently transferred from the Arabs to the Scandinavian merchants.

Many similar, so-called "regulated weights" were found by archeologists across the Scandinavian world, from Russia in the East to Ireland and England in the West [15]. In respect of their shape, most of the weights belong to two large groups: (1) the oblate-spheroid, sometimes referred to as 'spheres with flat poles', made of an iron core and copper-alloy mantle; (2) the copper-alloy cuboctahedron that is a main "personage" of our narrative.

There is a wide consensus in the archeological literature the weights of both groups reflect trading contacts with the Islamic Caliphate and most of such weights are considered to be imported from the East [15-16].

Obviously, for someone to be able to produce a cuboctahedral weight, this quite sophisticated shape must already be known. To the best of our knowledge, in archeological literature, these artifacts have never been considered with respect to the history of their geometric forms. In our opinion, the history of geometrical exploration of cuboctahedron and, more generally, of *Archimedean solids* can support the hypothesis on the Eastern origin of the cuboctahedral weight.

*Archimedean solids* are highly symmetric, semi-regular convex polyhedra composed of two or more types of regular polygons meeting in identical vertices (Figure 5, Table 1). The fact that *Archimedean solids* consist of at least 2 different types of polygons makes them distinct from the *regular polyhedra* or *Platonic solids* (cube, tetrahedron, octahedron, icosahedron, and dodecahedron). These polyhedra are called Archimedean because they were described by Archimedes (287 b. c. - 212 b. c.), even if we have only "second-hand" references to his writings on this topic from Heron of Alexandria (second half of the 1$^{st}$ century AD) and Pappus of Alexandria (290 - 350) [17-19].

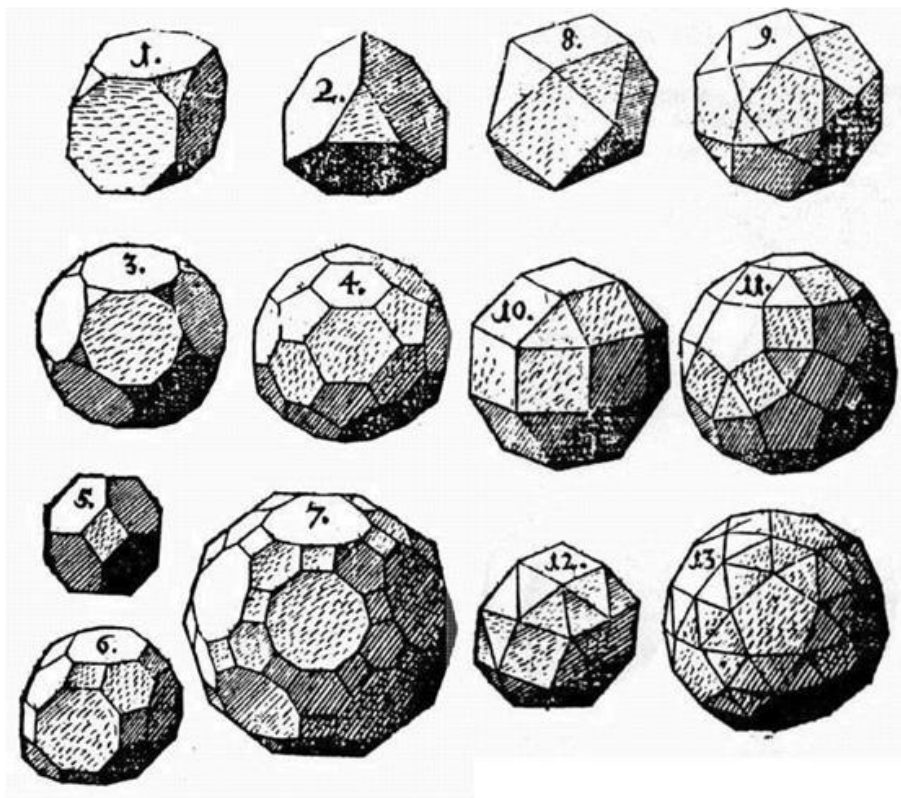

Figure 5. Images of Archimedean solids from Johannes Kepler's "Harmonice Mundi" [20]: polyhedron #8 is the cuboctahedron.

Table 1. Characteristics of the Archimedean polyhedra.

| Archimedean Solids | Number of | | | Rediscovered in Europe by |
|---|---|---|---|---|
| | Vertices | Faces[a] | Edges | |
| Truncated tetrahedron | 12 | $4_3 + 4_6 = 8$ | 18 | Piero della Francesca |
| Cuboctahedron | 12 | $8_3 + 6_4 = 14$ | 24 | Piero della Francesca |
| Truncated cube | 24 | $8_3 + 6_8 = 14$ | 36 | Piero della Francesca |
| Truncated octahedron | 24 | $6_4 + 8_6 = 12$ | 36 | Piero della Francesca |
| Rhombicuboctahedron | 24 | $8_3 + 6_4 + 12_4 = 26$ | 48 | Luca Pacioli |
| Snub cube | 24 | $24_3 + 8_3 + 6_4 = 38$ | 60 | Albrecht Dürer |
| Icosidodecahedron | 30 | $20_3 + 12_5 = 32$ | 60 | Luca Pacioli |
| Truncated cuboctahedron | 48 | $12_4 + 8_6 + 6_8 = 26$ | 72 | Albrecht Dürer |
| Truncated dodecahedron | 60 | $20_3 + 12_{10} = 32$ | 90 | Piero della Francesca |
| Truncated icosahedron | 60 | $12_5 + 20_6 = 32$ | 90 | Piero della Francesca |
| Rhombicosidodecahedron | 60 | $20_3 + 30_4 + 12_5 = 62$ | 120 | Luca Pacioli |
| Snub dodecahedron[b] | 60 | $60_3 + 20_3 + 12_5 = 92$ | 150 | Johannes Kepler |
| Truncated icosidodecahedron | 120 | $30_4 + 20_6 + 12_{10} = 62$ | 180 | Daniele Barbaro |

[a] The subindex indicates the order of the polygons that form the faces. For instance, for cuboctahedron, '$8_3 + 6_4$' stands for 8 triangles and 6 squares.

The knowledge of Platonic and Archimedean polyhedra was disseminated through the Arabic/Persian culture by translations made during the VIII-X centuries [19]. Construction of five *Archimedean solids* was described in the treatise of Persian mathematician and astronomer Abū al-Wafā Būzhjānī (940–998) "A Book on Those Geometric Constructions Which Are Necessary for a Craftsman", (كتاب في ما يحتاج إليه الصانع من الأعمال الهندسية) written after 990. It is important for our narrative that the cuboctahedron was among them [21].

Furthermore, the cuboctahedron was described at least a century earlier by Thābit ibn Qurra (826 or 836 – 901), another great mathematician and polymath from the Abbasid Caliphate, in his treatise "*On the Construction of a Solid Figure with Fourteen Faces Inscribed into a Given Sphere*". Most of the modern authors studied this manuscript by its German [22] and French [23] translations. Thabit ibn Qurra considered a spatial construction of a polyhedron bounded by six squares and eight equilateral triangles and illustrated this three-dimensional figure.

The geometry of cuboctahedral weights appears to have been transferred to certain types of silver jewelry in the second half of the IX century [15]. Cuboctahedrons were widely used in medieval Islamic architecture [24].

We suggest that cuboctahedron can be known earlier than other Archimedean solids since it can be produced by a mathematician's mind or by a craftsman's hands from a cube using a straightforward operation, so-called *rectification*: all cube's vertices should be truncated by cutting off a half of each edge of the cube adjacent to each vertex. It is probably not a coincidence that according to Heron, Archimedes ascribed the cuboctahedron to Plato (424/423 – 348 BC) [25].

In Europe, Archimedean polyhedra were rediscovered, described, and incorporated into the world of science and art only in the XV-XVII centuries (Table 1) [26].

Studying Nature, the Renaissance artists tried to find scientific ways of drawing it. Theories of perspective, proportions, light, and shade became a base for the new art. These theories were built on geometry and optics. They allowed an artist to create a three-dimensional space on a flat surface, saving an impression of relief of the objects. For some masters of the Renaissance, polyhedrons were just a convenient model for practicing the laws of perspective. Some of them were fascinated by their symmetry and laconic beauty. The others, following Plato, were attracted by philosophical and mystical symbols of polyhedral. It is worth noting that the treatise was also written by Abu al-Wafā' Būzjānī to provide craftsmen and artists with reasonable mathematical basics.

Five Archimedean polyhedra produced by truncation of Platonic solids were described by one of the greatest Renaissance artists and outstanding mathematician Piero della Francesca (1415-1492) between 1480 and 1490 in his treatise "Libbelus de quinque corporibus regularibus" ("Short book on the five regular solids") [27]. Piero's manuscript was lost for a long time but luckily, at the beginning of the XX century, its original was found. Nowadays it is in the Vatican Library. One of the co-authors (EAK) had a rare chance to keep it in his hands at the Vatican Library and was extremely happy to look at Piero's illustration of a truncated icosahedron. For most of his scientific career, this co-author has been investigating nanomaterials based on buckminsterfullerene, the molecule in the shape of a truncated icosahedron [28].

The main hero of our narrative, the cuboctahedron, was described by Piero in his other treatise "Trattato d'Abaco" ("Abacus Treatise") [29] written approximately at the same time as "Libbelus".

49A further input to the geometry of Archimedean solids was the book "De Divina Proportione" ("The Divine Proportion") written by the Franciscan friar and mathematician Luca Pacioli (1447-1517) between 1496 and 1498 printed in 1509 [30]. Pacioli described three polyhedra not mentioned in Piero's works.

Being Piero's countryman, Luca Pacioli started to study mathematics under his supervision. Some of the historians (the first was Giorgio Vasari who published a biography of Piero della Francesca in 1550 [31]) blame the author of "The Devine Proportion" for plagiarism of unpublished manuscripts of his teacher [3]. Indeed, Pacioli reproduced the descriptions of Archimedean solids from Piero's manuscripts without the reference. In any case, Pacioli's book strongly influenced the mathematics of the time. One of the main reasons for the book's success was the three-dimensional illustrations of polyhedra produced by Pacioli's friend Leonardo da Vinci [32].

Starting in 1525, i.e. three years before his death, Albrecht Dürer (1471-1528) hurried to share the secrets of perspective, that he had been acquiring all his life. He published two treatises, one of which, "Underweysung der Messung" ("Painter's Manual") [33] is a serious input in the theory of perspective and geometry of polyhedra. For instance, Dürer described two Archimedean solids unknown at this time.

Finally, two missing polyhedra were described by the Renaissance architect and writer (famous translator and commentator of Vitruvius) Daniele Barbaro (1514-1570) in the book "La pratica della perspettiva" (1568) [34] and the great Johannes Kepler (1571 - 1630) in his book of "Harmonice Mundi" ("Harmony of the Worlds") (1619) [20]. Kepler mathematically proved that there are only 13 polyhedra in the class of Archimedean solids [35-37], fully described each of them, and coined the names by which they are known today.

The rediscovery of Archimedean Solids lasted from 1480 till 1619. Though the cuboctahedron was one of the first such polyhedra rediscovered in Europe, one can conclude with certainty that by the time of the manufacture of the cuboctahedral weights discussed above, European mathematicians were not aware of such polyhedra, but Arab-Persian scientists and *craftsmen* were.

**References and commentaries**

1. J J. Gray and P.R. Cromwell. Kepler's Work on Polyhedra, *Mathematical Intelligencer*, 17 (1995), 23–33.
2. M. A. Peterson. The geometry of Piero della Francesca, *Mathematical Intelligencer*, 19 (1997), 33-40.
3. K. Williams. Plagiary in the Renaissance, *Mathematical Intelligencer*, 24, (2002), 45–57.
4. D. Huylebrouck. Lost in Enumeration: Leonardo da Vinci's Slip-Ups in Arithmetic and Mechanics, *Mathematical Intelligencer*, 34 (2012), 15-20.
5. E. A. Katz and B. Y. Jin. Fullerenes, Polyhedra, and Chinese Guardian Lions, *Mathematical Intelligencer*, 38 (2016), 61-68.
6. Cultural layer is a key concept in archaeology. It refers to remnants of human settlement (details of housing structures and farming, tools, ceramics, sacrifices, or other indicators of rituals) often buried by centuries or millennia of sediment and discovered during archaeological excavations.
7. Н. Б. Черных. Дендрохронология древнейших горизонтов Старой Ладоги (по материалам раскопки Земляного городище. *Средневековая Ладога*, Ленинград, Л., 1985, 76-80.
8. Е. А. Рябинин Н. Б. Черных. Стратиграфия, застройка и хронология нижнего слоя Староладожского Земляного городища в свете новых исследований, *Советская археология* (1988), №1, 72-100. https://www.archaeolog.ru/media/books_sov_archaeology/1988_book01.pdf
9. *Полное собрание русских летописей* Т. II. 2-е изд. Санкт Петербург, 1908.
10. *Повесть временных лет.* Вита Нова, Санкт Петербург, 2012. https://imwerden.de/pdf/povest_vremennykh_let_2012.pdf
11. M. Dimnik, (1994). *The Dynasty of Chernigov 1054–1146*. Pontifical Institute of Mediaeval Studies, 1994.
12. H. G. Lunt. What the Rus' Primary Chronicle Tells Us about the Origin of the Slavs and of Slavic Writing. *Harvard Ukrainian Studies,* 19, (1995) 335–357.
13. Е. А. Рябинин. Начальный этап стеклоделия в Балтийском регионе (по материалам исследований Ладоги VIII—X вв.). *Дивинец Староладожский*. СПб., 1997, 43-49.

14. А. Ю. Чернов. Энеида Рёрика Ютланского (2022).

    https://imwerden.de/pdf/chernov_eneida_rerika_yutlanskogo_2022__txt.pdf

15. J. Kershaw. Metrology and beyond: new approaches to Viking-Age regulated weights. Weight and Value, 1, (2019) 121-132.

16. U. Pedersen, In: *Means of Exchange: Dealing with Silver in the Viking Age*. Kaupang Excavation Project Publication, Series 2: Norske Oldfunn XXIII, Oslo, 2008, 119-178.

17. The manuscript of Archimedes is generally thought to have been lost on the famous conflagration of the Ancient Library of Alexandria. There was a reference to this manuscript and description of Archimedean solids in the book (collection of 8 books) of one of the last great Greek mathematicians of antiquity Pappus of Alexandria (290 – 350) "Synagoge" "or "Collection" ("Synagoge" is "Collection" in Greek, συναγωγή). Alexander Jones wrote [18]: "The only ancient source for these solids [the Archimedean solids] is Pappus Book 5, chapters 34-37, together with a marginal note describing the construction of some of them". Wide knowledge of this book had to wait until the end of XVI century when its Latin translation was published.

18. A. Jones, *Pappus of Alexandria: Book 7 of the Collection, 2 volumes*. Springer, 1986, volume I, 50.

19. P. R. Cromwell, *Polyhedra*, Cambridge University Press, Cambridge, 1999.

20. J. Kepler, *Harmonice Mundi*, Linz, 1619.

21. R. Sarhangi. Illustrating Abu al-Wafā' Būzjānī: Flat Images, Spherical Constructions, *Iranian Studies*, 41, (2008) 511-523.

22. E. Bessel-Hagen and O. Spies. Ṯābit b. Qurra's Abhandlung über einen halbregelmäßigen Vierzehnflächner'. *Quellen und Studien zur Geschichte der Mathematik, Astronomie und Physik*. Abteilung B: Studien 2 (1933), 186–198.

23. K. Asselah. Thâbit ibn Qurra: Construction d'un polyèdre semi-régulier à quatorze faces, 8 trianges équilatéraux et 6 carrés / Texte et traduction: Construction d'une figure solide à quatorze faces. In: *Thabit ibn Qurra: Science and Philosophy in Ninth-Century Baghdad*, eds. Roshdi Rashed, 317-323/ā324-331 (2009). Berlin: Walter de Gruyter.

24. H. Hisarligil, B. B. Hisarligil. The Geometry of Cuboctahedra in Medieval Art in Anatolia. *Nexus Network Journal,* 20 (2018), 125–152.

25. H.S.M. Coxeter, *Regular Polytopes*. 1973, New York: Dover Publications.

26. J.V. Field. Rediscovering the Archimedean Polyhedra: Piero della Francesca, Luca Pacioli, Leonardo da Vinci, Albrecht Durer, Daniele Barbaro, and Johannes Kepler, *Archive for History of Exact Sciences*, 50, (1997) 241-283.
27. Piero della Francesca. *Libbelus de quinque corporibus regularibus*. Vatican library, ttps://digi.vatlib.it/mss/detail/Urb.lat.632 .
28. E.A. Katz. Bridges between mathematics, natural sciences, architecture and art: case of fullerenes. In: *Proc. of the 1st International Conference "Art, Science and Technology: Interaction between Three Cultures"*, Domus Argenia Publisher, Milano, Italy, 60-71, 2012.
29. P. della Francesca, *Trattato d'Abaco*, ed. G. Arrighi, Pisa: Domus Galilaeana (1970).
30. L. Pacioli, *De Divine Proportione*, 1509, (Ambrosiana fascimile reproduction, 1956; Silvana fascimile reproduction, 1982).
31. G. Vasari, *The lives of the artists*, Oxford University Press, 1991.
32. D. Huylebrouck, Observations about Leonardo's Drawings for Luca Pacioli, *BSHM Bulletin: Journal of the British Society for the History of Mathematics*, 30, (2014) 102-112.
33. A. Dürer, *Underweysung der Messung*, 1525 (translated to English as "Painter's Manual", Abaris reprint, 1977).
34. D. Barbaro, *La pratica della perspettiva*, Venice, 1568.
35. A new wave of interest in the polyhedra in XX century is related to the role that symmetry considerations play in modern science. In particular, an "additional" semiregular body, evidently overlooked by Archimedes or Kepler, was discovered independently by several scientists in different countries in the XX century (see, for example [32]). This polyhedron called pseudo-rhombicuboctahedron can be obtained from the rhombicuboctahedron (number 10 in Figure 2) by turning one octagonal bowl by 45 ° with respect to the entire body. Such new body satisfies to the definition of Archimedean solids given at the beginning of this essay. It can be considered as the 14th Archimedean body. However, M. J. Wenninger, the author of "Polyhedron Models" [33] does not believe that it is the Archimedean body because square faces with cubic and rhombic origin are mixed "confused" in it.
36. B. Grünbaum. An enduring error, *Elemente der Mathematik*, 64 (2009) 89–101. Reprinted in Pitici, Mircea, ed. (2011), The Best Writing on Mathematics 2010, Princeton University Press, pp. 18–31.
**37.** M. J. Wenninger, *Polyhedron Models*, 1971, Cambridge University Press.